\documentclass[a4paper,12pt]{amsart}
  \usepackage{amssymb,amsthm}
\usepackage{graphicx,color}        % standard LaTeX graphics tool
\usepackage{url} % <-- Para p\'aginas web o similar: \url{...}

\newtheorem{theorem}{Theorem}[section]
\newtheorem{corollary}[theorem]{Corollary}
\newtheorem{proposition}[theorem]{Proposition}

\theoremstyle{definition}
\newtheorem{definition}[theorem]{Definition}

\newtheorem{remark}[theorem]{Remark}

\def\r{\mathbb R}
\def\s{\mathbb S}
 \def\h{\mathbb H}
 \def\n{\mathbf n}

\begin{document}

\title{The hanging chain problem in the sphere and in the hyperbolic plane}
\author{Rafael L\'opez}
\address{ Departamento de Geometr\'{\i}a y Topolog\'{\i}a\\  Universidad de Granada.  18071 Granada, Spain}
\email{rcamino@ugr.es}
%
% Use the package "url.sty" to avoid
% problems with special characters
% used in your e-mail or web address
%

\begin{abstract} In this paper, the notion of the catenary curve in the sphere and in the hyperbolic plane is introduced. In both spaces, a catenary is  defined as the shape of a hanging chain when its  potential energy is determined by the distance to a given geodesic of the space. Several characterizations of the catenary are established in terms of the curvature of the curve  and of the angle that its unit normal makes with a vector field of the ambient space. Furthermore, in the hyperbolic plane, we extend the concept of catenary substituting the reference geodesic by a horocycle or the hyperbolic  distance by the horocycle distance.
\end{abstract}

\maketitle

\section{Introduction and objectives}
%%%%%%%%%%%%%%%%%

The shape that adopts a hanging chain under its own weight when suspended from its endpoints attracted the interest of scientistics from times of Galileo and da Vinci. Galileo   believed that the parabola was the shape of the chain  but his argument was wrong. The solution curve is not so simple as a parabola (a quadratic polynomial function) but the catenary 
\begin{equation}\label{cat}
y(x)=\frac{1}{a}\cosh(ax+b)-\lambda,\quad a,b,\lambda\in\r, a>0,
\end{equation}
which is a curve involving transcendental functions as the exponential.   The derivation of the solution was an independent work of  R. Hooke, J. Bernouilli, G. Leibniz and C. Huygens among others. See \cite{beh,co} for an account of the history of the catenary. The catenary is a classical curve and one of the first examples, together the brachistochrone,  that illustrates the power of the calculus of variations \cite{gh}. Related with the catenary, and also using calculus of variations, Euler proved that the catenary is the generating curve of the surface of revolution with minimum area and spanning to coaxial circles \cite{eu}.  To be precise, if the catenary \eqref{cat} rotates around the $x$-axis, the resulting surface of revolution is minimal if and only if $\lambda=0$. For other mathematical properties of the catenary,  see \cite{ch,cd,kkp,mc2,pa}.

The catenary appears related with different topics in science, specially in engineering and architecture. For example, it is the model of an arch where the only force acting on the arch is its weight \cite{he,pot}. This suggests its utilization in the construction of arches and roofs of corridors, such as the Spanish architect A. Gaud\'{\i} used in many of its constructions, as for example,   the Colegio Teresiano and La Pedrera (Barcelona). In this sense, it is very nice to read the two articles of R. Osserman about the shape of the Gateway Arch in St. Louis, Missouri, connecting the shape of the Gateway Arch and the catenary \cite{os,os2}. 

Many extensions of the hanging chain problem have been investigated and the literature is huge. Without to give a complete list of references, we point out some of the modifications in the classical problem.  For example, one can assume that: the density of the chain changes along its length \cite{fa,kk,ok};   the force vector field is radial  \cite{dh};  the chain is made of an elastic material  \cite{bow,ir2,ir,irsin,rl};  the chain is  subjected   under the effect of the surface tension of a soap film   adhered to the chain \cite{bmm1,bmm2};    the chain is suspended from a vertical line and rotates around this axis  \cite{ap,ms,ne};  the two ends of the hanging  chain move  with stretching the chain along a path \cite{kaj};   there are loads on  the chain  which pulling down on its lowest point  \cite{za}.

In this paper we will investigate the generalization of the concept of  catenary in the unit sphere $\s^2$ and in  the hyperbolic plane $\h^2$. From the mathematical viewpoint, it is natural to ask for the extension of the hanging chain problem to other spaces. It is clear that the sphere and the hyperbolic plane are the first spaces to study because they are the models of the elliptic geometry and the hyperbolic geometry, respectively. However, it is surprising that this theme has not been considered in the literature, being the catenary   well-known for centuries, as well as the sphere and the hyperbolic plane are classical models in geometry. 

The purpose of this article is to give an approach to the generalization of the hanging chain problem in these two spaces. More specifically,
 to formulate a suitable problem that can be adopted as an extension of the Euclidean catenary. Once the concept of catenary is defined   in these spaces  as a critical point of a potential energy functional,   different characterizations of the solution curve  in terms of its curvature will be obtained. Finally, and if possible, we ask if the curves obtained as solutions of the hanging problem are the generating curves of minimal surfaces in the $3$-dimensional sphere $\s^3$ and hyperbolic space $\h^3$. This would extend the Euler's result to these spaces. 
 
Our motivation to generalize the notion of the catenary in the sphere and in the hyperbolic plane has its origin in the Euclidean catenary \eqref{cat}. We will recall the hanging chain problem in the Euclidean plane $\r^2$, pointing out which are the ingredients in its formulation. These can give us the clues to proceed when the ambient space is the sphere and  the hyperbolic plane.

 Let $\r^2$ denote the Euclidean plane where $(x,y)$ stand for the Cartesian coordinates and let $\langle,\rangle$ be the Euclidean metric.   The hanging chain problem consists in finding the shape of an inextensible chain with uniform linear mass density and  suspended from two fixed endpoints. Suppose that the chain of mass $m$ is idealized as a curve  $y\colon [a,b]\to\r^2$, $y=y(x)$.  The gravitational acceleration $g$ is constant over the chain. The $x$-axis is taken as the level of zero potential energy. The gravitational potential energy of an infinitesimal element $ds$ of the chain at $(x,y)$ is $gy\, dm=\sigma g y\, ds$, where $\sigma$ is the density per unit length. Since $ds=\sqrt{1+y'(x)^2}\, dx$, the total potential energy of the chain is
\begin{equation}\label{eq0}
 \int_a^b\sigma g y(x)\sqrt{1+y'(x)^2}\ dx.
\end{equation}
  We are assuming  in \eqref{eq0} that $y(x)>0$ for all $x\in[a,b]$.  The hanging chain problem reduces to find a curve $y=y(x)$ that minimizes this energy among all curves with the same ends and the same length. The latter hypothesis is due to the inextensibility of the chain and the absence of elastic forces. Simplifying  the constant $\sigma g$ by $1$, the energy functional to minimize is
\begin{equation}\label{ee}
\mathcal{E}[y]= \int_a^b y\sqrt{1+y'^2}\, dx +\lambda \int_a^b \sqrt{1+y'^2}\, dx .
\end{equation}
The second term of $\mathcal{E}[y]$ is a Lagrange multiplier because all curves have the same length. Consequently, the solution $y(x)$ is a critical point of the energy  $\mathcal{E}$.  Using standard  arguments of calculus of variations, the Euler-Lagrange equation of \eqref{ee} is 
\begin{equation}\label{catenary}
\frac{y''}{(1+y'^2)^{3/2}}=\frac{1}{(y+\lambda)\sqrt{1+y'^2}}.
\end{equation}
The solution of \eqref{catenary} is  the    catenary  \eqref{cat}.  Note that the left-hand side of \eqref{catenary} is  the curvature $\kappa_e$ of the plane curve $y=y(x)$. The right-hand side of \eqref{catenary}   has the following  geometric interpretation.   Consider the vector field $\partial_y$ which is the (constant) gravitational field in $\r^2$.  Since the unit normal vector ${\bf n}$ of the curve $y(x)$  is 
${\bf n}=(-y',1)/\sqrt{1+y'^2}$, then $\langle {\bf n},\partial_y\rangle=1/\sqrt{1+y'^2}$ and  equation \eqref{catenary} can be expressed as
\begin{equation}\label{catenary2}
\kappa_e=\frac{\langle{\bf n},\partial_y\rangle}{y+\lambda}.
\end{equation}
Equation \eqref{catenary2}  shows that the hanging chain problem is equivalent to a coordinate-free prescribed curvature problem.

Motivated by the above description of the problem,  the main ingredients are the following.  The first aspect concerns to the existence of  a  reference line which is prescribed in the problem. In the above arguments, this line is the $x$-axis of $\r^2$ and at this level, the potential is $0$. A second aspect is the existence of a potential energy which depends on the position with respect to the reference line. In the case of the catenary, this potential is due to the gravity. It is to this potential that we want to calculate a minimum energy, or more exactly, a critical point. Finally, in the variational problem, all curves of the variation have prescribed endpoints and the same length. In particular, it is necessary to add a Lagrange multiplier to the potential energy that we want to minimize. 
Based on the above discussion, we will formulate the hanging chain problem in the sphere $\s^2$ and in the hyperbolic plane $\h^2$.  

The objectives of this paper can be divided into three specific items:  
\begin{enumerate}
\item[(T1)] State the analogous hanging chain problem in $\s^2$ and in $\h^2$ and find the corresponding Euler-Lagrange equation. 
\item[(T2)] Obtain an analogous formulation of the prescribed curvature equation \eqref{catenary2} in terms of a vector field  that represents the   `gravitational vector field'.
\item[(T3)]   Rotate the catenary in $\s^3$ and $\h^3$ and determine any unique properties of the mean curvature of the resulting surface.  \end{enumerate}

 The critical points of the potential energy will be also called  catenaries. Catenaries in $\s^2$ will be discussed in Section \ref{sec2} where two potential energies are used, first with the distance to a geodesic of $\s^2$ and second with the distance to a plane of $\r^3$. In the hyperbolic plane $\h^2$, the reference lines will be geodesics as well as horocycles.  This work  is   carried out in Section \ref{sec3}.

%%%%%%%%%%%%%%%%%%%%%%%%%%%%%%%
\section{The hanging chain problem in the sphere}\label{sec2}
%%%%%%%%%%%%%%%%%%%%%%%%%%%%%%%%%%%%%%%%%%%%%%%%%%%%%%%

In this section we will consider the hanging chain problem in the unit sphere $\s^2$. We first state the problem, and then give different characterizations of its solution. 

\subsection{Spherical catenaries: definition}

Consider the unit sphere $\s^2=\{(x,y,z)\in\r^3:x^2+y^2+z^2=1\}$. Let $\Psi$ be  the standard parametrization of $\s^2$ given by 
\begin{equation}\label{para}
\Psi(u,v)=(\cos{u}\cos{v},\cos{u}\sin{v},\sin{u}),\quad u\in [-\frac{\pi}{2},\frac{\pi}{2}],v\in\r.
\end{equation}
In the hanging chain problem in $\s^2$, consider a geodesic of $\s^2$ as the reference line to calculate the potential energy of a chain contained in $\s^2$. Let us fix the great circle $P=\{(x,y,z)\in\s^2:z=0\}$ as the reference geodesic.  Notice that on $\s^2$ we have not a concept of gravity $g$ due to the curvature of $\s^2$. 

At this point, we assign to each point of $\s^2$ a potential which measures its distance to the reference line $P$, being $P$ the level of zero potential. This distance is realized along all geodesics (meridians) orthogonal to $P$. The distance $d$ of a point $(x,y,z)=\Psi(u,v)$ to the geodesic $P$ is  $d=|\arcsin(z)|=|u|$. We also consider the unit vector field $X\in\mathfrak{X}(\s^2)$ (except at the north and south poles) which is tangent to all these geodesics. This vector field is the gradient $\nabla d$ of the distance function, which it is $\frac{\partial\Psi}{\partial u}=\Psi_u$. The   vector field $X$ can be expressed in terms of the canonical vector fields $\{\partial_x,\partial_y,\partial_z\}$ of $\r^3$ as  
\begin{equation}\label{x1}
X(\Psi(u,v))=-\sin{u}\cos{v}\, \partial_x -\sin{u}\sin{v}\, \partial_y+\cos{u}\, \partial_z.
\end{equation}
As in the Euclidean case, we will consider curves of $\s^2$ that do not intersect $P$, hence $u\not=0$. The geodesic $P$ separates $\s^2$ in two domains, namely,  the half-spheres  $\s^2_+=\{(x,y,z)\in\s^2:z>0\}$ and $\s^2_{-}=\{(x,y,z)\in\s^2:z<0\}$. Without loss of generality,  we will assume that all curves will be contained in the upper half-sphere $\s^2_+$. Let $\gamma  \colon [a,b]\to\s^2_{+}$ be a regular curve.  Let us write $\gamma(t)=\Psi(u(t),v(t))$, $t\in [a,b]$,  with the condition  $u(t)\in (0,\pi/2]$ because $\gamma(t)\in\s^2_{+}$. The arc-length element of $\gamma$  is  $\sqrt{u'^2+v'^2(\cos{u})^2}\, dt$. Consequently,  the potential energy of $\gamma$ is  
\begin{equation}\label{esg}
\mathcal{E}_S[\gamma]=\int_a^b (u+\lambda)\, |\gamma'(t)|\, dt=\int_a^b (u+\lambda)\sqrt{u'^2+v'^2(\cos{u})^2}\, dt, 
\end{equation}
where $\lambda\in\r$. The second integral of $\mathcal{E}_S[\gamma]$ is a Lagrange multiplier because in the variational problem   all curves   have the same length. 

\begin{definition} \label{ds1}
A critical point of $\mathcal{E}_S$ is called a {\it spherical catenary}.
\end{definition}

Before to find the critical points of  $\mathcal{E}_S$,  we will obtain a suitable expression for the curvature  of a curve in $\s^2$. Here the curvature  is understood  to be the geodesic curvature $\kappa_s$ of $\gamma$ in $\s^2$. The sign of $\kappa_s$ depends on the orientation of $\s^2$ which will be $N(p)=-p$, $p\in\s^2$. In such a case, 
\begin{equation}\label{ks}
\begin{aligned}
\kappa_s&=\frac{\langle \gamma'',\gamma'\times N(\gamma)\rangle}{|\gamma'|^3}=\frac{\langle \gamma,\gamma'\times \gamma''\rangle}{|\gamma'|^3}\\
&=\frac{1}{|\gamma'|^3}\left(v'(2u'^2\sin{u}+v'^2(\cos{u})^2\sin{u})-\cos{u}(u'v''-u''v')\right).
\end{aligned}
\end{equation}

Related to the objective T1, we have the following result.

\begin{theorem}\label{ts1}
 Let $\gamma(t)=\Psi(u(t),v(t))$ be a regular curve in $\s^2_+$. Then $\gamma$ is a spherical catenary  if and only if its   curvature $\kappa_s$  satisfies
\begin{equation}\label{s2}
 \kappa_s=  \frac{v'\cos{u}}{(u+\lambda)|\gamma'|}.
\end{equation}
\end{theorem}

\begin{proof} 
We calculate the Euler-Lagrange equation of the energy \eqref{esg}. The  Lagrangian of $\mathcal{E}_S$   is 
$J[u,v,u',v']=(u+\lambda)\sqrt{u'^2+v'^2(\cos{u})^2}$. Since the same computations will be done later in a similar context (see Proposition \eqref{pr1}), we assume a more general case of Lagrangian of type 
\begin{equation}\label{lagrange}
J[u,v,u',v']=f(u)\sqrt{u'^2+v'^2(\cos u)^2}.
\end{equation}
A curve $\gamma$ is a critical point   if and only if $\gamma$ satisfies  
\begin{equation}\label{sel}
\frac{\partial J}{\partial u}-\frac{d}{dt} \left(\frac{\partial J}{\partial u'}\right)=0\quad\mbox{and}\quad 
 \frac{\partial J}{\partial v}-\frac{d}{dt} \left(\frac{\partial J}{\partial v'}\right)=0.
\end{equation}
After some  computations, equations \eqref{sel} are, respectively, 
\begin{equation}\label{eqs0}
\begin{aligned}
\frac{v'\cos{u}}{|\gamma'|}\left(f'v'\cos{u}-fu'\sin{u}\right)-f\frac{d}{dt}\left(\frac{u'}{|\gamma'|}\right)&=0,\\
 \frac{f'u'v'(\cos{u})^2}{|\gamma'|}+f\frac{d}{dt}\left(\frac{v'(\cos{u})^2}{|\gamma'|}\right)&=0.
 \end{aligned}
 \end{equation}
Equations \eqref{eqs0} can be written in terms of $\kappa_s$ as follows. Using  \eqref{ks}, we obtain
\begin{equation}\label{eqs1}
\begin{aligned}
\frac{d}{dt}\left(\frac{u'}{|\gamma'|}\right)&=\frac{v'\cos{u}}{|\gamma'|^3}\left(v'u'^2\sin{u}-\cos{u}(u'v''-u''v')\right)\\
&=v'\cos{u}\left(\kappa_s-\frac{v'\sin{u}}{|\gamma'|}\right),
\end{aligned}
\end{equation}
and 
\begin{equation}\label{eqs2}
\begin{aligned}
\frac{d}{dt}\left(\frac{v'(\cos{u})^2}{|\gamma'|}\right)&=-\frac{u'\cos{u}}{|\gamma'|^3}\left(2u'^2v'\sin{u}+v'^3\sin{u}(\cos{u})^2-\cos{u}(u'v''-u''v')\right)\\
&=- u'  \kappa_s \cos{u}.
\end{aligned}
\end{equation}
From \eqref{eqs1} and \eqref{eqs2}, the Euler-Lagrange equations \eqref{eqs0} can be expressed as
\begin{equation}\label{uuvv}
\begin{aligned}
v'\cos{u}\left(\frac{f'v'\cos{u}}{f|\gamma'|}-\kappa_s\right)&=0,\\
u'\cos{u}\left(\frac{f'v'\cos{u}}{f|\gamma'|}-\kappa_s\right)&=0.
\end{aligned}
\end{equation}
Since  $\gamma$ is a regular curve,  $u'$ and $v'$ cannot be simultaneously zero.   Together with $\cos{u}\not=0$ and \eqref{uuvv}, we deduce
\begin{equation}\label{silva}
\kappa_s=\frac{f'v'\cos{u}}{f|\gamma'|}.
\end{equation}
In the particular case that $f(u)=u+\lambda$, then \eqref{silva} is just \eqref{s2}.  
\end{proof}

Equation \eqref{s2} is second order, but a first integration is possible because the Lagrangian $J[u,v,u',v']$ does not depend on the function $v$. Indeed, there is a constant $c$ such that 
\begin{equation}\label{si}
\frac{\partial J}{\partial v'}= (u+\lambda) \frac{v'(\cos{u})^2}{\sqrt{u'^2+v'^2(\cos{u})^2}}=c.
\end{equation}
If $c=0$, then $v'=0$ and $\gamma$ is a meridian of $\s^2$. Thus, if $c\not=0$, then $\gamma$ is never tangent to a meridian. In particular, $\gamma$ does not across the north pole. Without loss of generality, we can assume  that  $u$ is a function of $v$, $u=u(v)$. Then  $\gamma(v)=\Psi(u(v),v)$ and  \eqref{si} can be rewriten as  
\begin{equation*} 
(u+\lambda) \frac{ (\cos{u})^2}{\sqrt{u'^2+ (\cos{u})^2}}=c.
\end{equation*}
Hence, an expression for $u=u(v)$ is deduced, obtaining 
\begin{equation}\label{deduced}
\int^u\frac{du}{\cos{u}\sqrt{(u+\lambda)^{2}(\cos{u})^2-c^2}}=\frac{v}{c}.
\end{equation}
This integral yields the following corollary.

\begin{corollary}
Let $\gamma(v)=\Psi(u(v),v)$ be a regular curve in $\s^2_+$. Then $\gamma$ is a spherical catenary  if and only if   $u=u(v)$ satisfies \eqref{deduced} for some constant $c\in\r$. 
\end{corollary}

\begin{remark} Identity \eqref{si} is a type of Clairaut relation for spherical catenaries. Since the Clairaut relation on $\s^2$ holds for geodesics, we need to eliminate the Lagrange constraint due to the length in the initial formulation of the variational problem. Thus, take $\lambda=0$ in \eqref{si}. Under this assumption, the angle $\Theta$ that $\gamma$ makes with the parallel $v\mapsto \Psi (u,v)$ is 
$$\cos\Theta=\frac{\langle\gamma',\Psi_v\rangle}{|\gamma'||\Psi_v|}=\frac{v'\cos{u}}{\sqrt{u'^2+v'^2(\cos{u})^2}}.$$
Then identity \eqref{si} can be expressed as
\begin{equation}\label{clairaut}
u \cos{u}\cos\Theta=c.
\end{equation}
The classical Clairaut relation establishes $\cos{u}\cos\Theta=c$ \cite[p. 257]{carmo}. In the case of spherical catenaries, identity \eqref{clairaut} asserts that the radius of the parallel multiplied by $u$ and by the cosine of the intersection angle with each parallel is constant.
\end{remark}

\begin{figure}[hbtp]
\centering
\includegraphics[width=.3\textwidth]{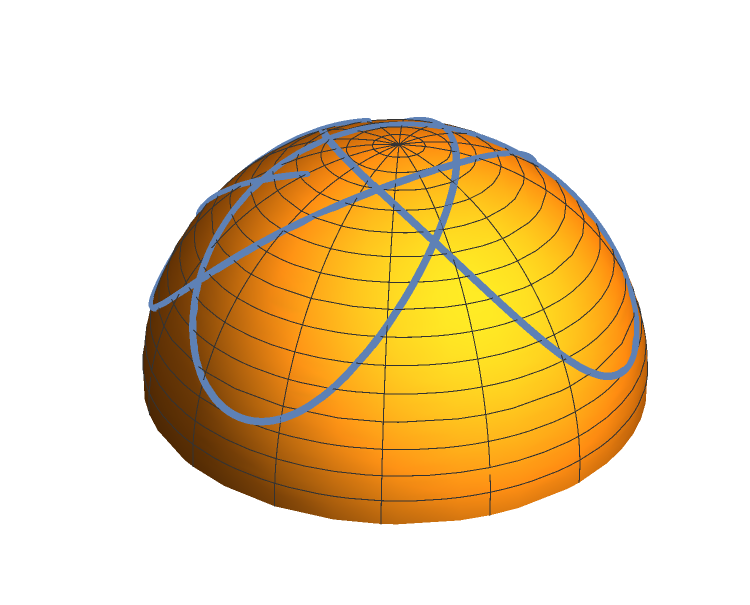} \includegraphics[width=.3\textwidth]{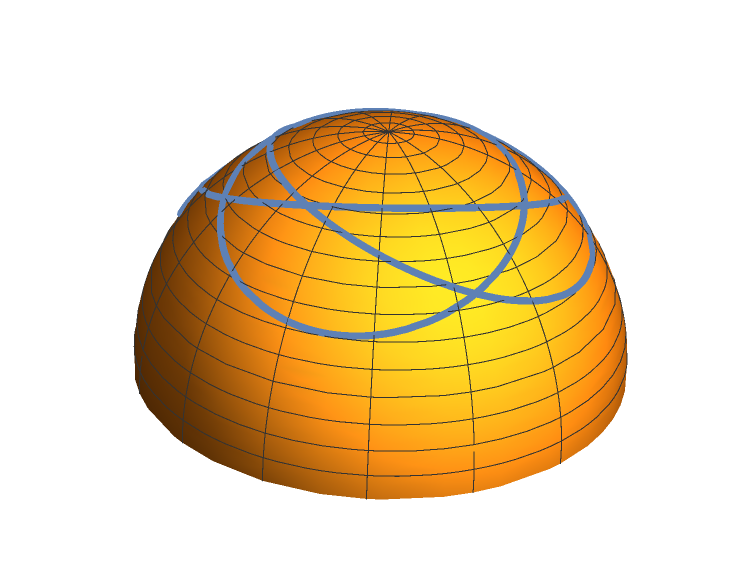} \includegraphics[width=.3\textwidth]{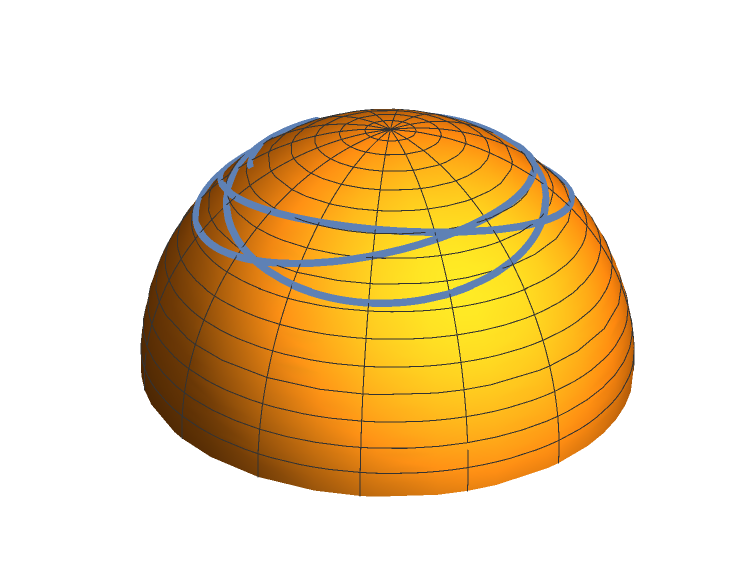}
\caption{Spherical catenaries in $\s^2_{+}$ (top view from the north pole) for different initial conditions in the ODE system \eqref{h00}. The initial conditions are $v(0)=0$, $\theta(0)=\pi/2$ and different values for $u(0)$: $u(0)=0.3$ (left), $u(0)=0.6$ (middle) and $u(0)=1$ (right).}\label{fig1}.
\end{figure}

Going back, equation \eqref{s2} can be viewed as a prescribed equation of a curve $\gamma$ in the sphere. Since  $\gamma$ is the image of the plane curve $\beta(t)=(u(t),v(t))$ under the parametrization \eqref{para}, equation \eqref{s2} can be reformulated in terms of the curvature $\kappa_\beta$ of $\beta$. The curvature $\kappa_\beta$ of $\beta$ is 
\begin{equation}\label{k5}
\kappa_\beta(t)=\frac{u'v''-u'' v'}{(u'^2+v'^2)^{3/2}}.
\end{equation}
Inserting in \eqref{ks}, we have 
$$\kappa_s=\frac{1}{|\gamma'|^3}\left( 2u'^2v'\sin{u}+v'^3(\cos{u})^2\sin{u}-\kappa_\beta\cos{u}(u'^2+v'^2)^{3/2}\right).$$
Thus  equation \eqref{s2} can be written as
\begin{equation}\label{kapabeta}
\kappa_\beta=\frac{2u'^2v'\sin{u}+v'^3(\cos{u})^2\sin{u}-\frac{v'\cos{u}}{u+\lambda}\left(u'^2+v'^2(\cos{u})^2\right)}{\cos{u}(u'^2+v'^2)^{3/2}}.
\end{equation}
 This expression for $\kappa_\beta$ allows to illustrate some spherical catenaries in Figure \ref{fig1}. These plots have been made with the {\it Mathematica}  software  (\cite{wo}).   We briefly explain  the method to obtain these figures.  Suppose $\gamma(t)=\Psi(u(t),v(t))$ and that the curve $\beta(t)=(u(t),v(t))$ is parametrized by arc-length. Then $u'(t)^2+v'(t)^2=1$. So we can write $u'(t)=\cos\theta(t)$ and $v'(t)=\sin\theta(t)$ for some function $\theta=\theta(t)$. According to \eqref{kapabeta}, the functions $u(t)$ and $v(t)$ satisfy the ODE system, 
\begin{equation}\label{h00}
\left\{ \begin{aligned}
u'(t)&=\cos\theta(t)\\
v'(t)&=\sin\theta(t)\\
\theta'(t)&=\sin\theta\left( 2(\cos\theta)^2 \tan{u}+(\sin\theta)^2\cos{u}\sin{u}-\frac{1-(\sin\theta)^2(\sin{u})^2}{u+\lambda} \right).
\end{aligned}
\right.
\end{equation}
Recall that the variation of the angle function $\theta(t)$ coincides with the curvature $\kappa_\beta$ of $\beta$.  Given initial conditions $u(0)=u_0$, $v(0)=v_0$ and $\theta(0)=\theta_0$, {\it Mathematica} solves numerically the ODE \eqref{h00} and then the software graphically represents the solution $\gamma(t)= \Psi(u(t),v(t))$.  
 
%%%%%%%%
\subsection{Spherical catenaries: characterizations}

 In this subsection, we give a geometric interpretation of \eqref{s2} completing the objective T2.   

\begin{theorem}\label{ts2}
Let $\gamma$ be a regular curve in $\s^2_{+}$. Then $\gamma$ is a spherical catenary if and only if its geodesic curvature $\kappa_s$ satisfies
 \begin{equation}\label{s3}
 \kappa_s =  \frac{\langle{\bf n},X\rangle}{d+\lambda},
\end{equation}
where ${\bf n}$ is the unit normal vector of the curve $\gamma$ (as a tangent vector on $\s^2$ orthogonal to $\gamma'$), and $d$ is the distance to $P$.
\end{theorem}

\begin{proof}
  Since the unit normal vector to $\s^2$ along $\gamma$ is  $N=-\gamma$,   the  vector ${\bf n}$ is 
\begin{equation}\label{nn}
{\bf n}=\frac{ \gamma'\times N(\gamma)}{|\gamma'|}=\frac{\gamma\times\gamma'}{|\gamma'|}.
\end{equation}
From the definition of $X$,   $\langle {\bf n},X(\gamma)\rangle= v'\cos{u}/|\gamma'|$. Equation \eqref{s3}  follows from \eqref{s2} because  the function $u$ in \eqref{s2} is just the distance $d$ to $P$.
 \end{proof}

Notice that  equation \eqref{s3} is the analogue of   \eqref{catenary2}  for spherical catenaries of $\s^2$.  
 
The last part of this section addresses the third objective T3. Minimal rotational surfaces of $\s^3$ have been studied in the literature. However,   there is no   a known geometric property of the generating curves of these surfaces. More exactly, we ask if the generating curves can be viewed as solutions of a hanging chain problem in $\s^3$. 

A surface of revolution in the $3$-dimensional sphere $\s^3$ will be constructed by rotating a spherical catenary  around the geodesic $P$. Here, as in the Euclidean catenary,  we will assume $\lambda=0$ for the Lagrange multiplier. Let $\s^3=\{(x_1,x_2,x_3,x_4)\in\r^4:x_1^2+x_2^2+x_3^2+x_4^2=1\}$ and   $\s^2\hookrightarrow\s^3$ be the natural inclusion defined by $(x_1,x_2,x_3)\mapsto (x_1,x_2,x_3,0)$. This embedding identifies $\s^2$ with $\s^2\times\{0\}\subset\s^3$.  Let $\gamma\colon I\to \s^2\subset\s^3$ be a curve contained in $\s_+^2$.  Denote by $S_\gamma$ the surface of revolution in $\s^3$ obtained by rotating $\gamma$ with respect to $P\subset\s^2\times\{0\}$.  The one-parameter group of rotations whose axis is $P$ is $\mathcal{G}=\{\mathcal{R}_s:s\in\r\}$, where
$$\mathcal{R}_s=\left(\begin{array}{cccc}1&0&0&0\\ 0&1&0&0\\ 0&0&\cos{s}&-\sin{s}\\ 0&0&\sin{s}&\cos{s}\end{array}\right).$$

 In order to simplify the computations, we can assume without loss of generality that   $\gamma$ is parametrized by $\gamma(t)=\Psi(u(t),t)$. Then the parametrization of $S_\gamma$  is 
\begin{equation}\label{para1}
\begin{aligned}
\Phi(t,s)&=\mathcal{R}_s\cdot\gamma(t)\\
&=\left(\cos{(u(t))}\cos t,\cos{(u(t))}\sin{t},\sin{(u(t))}\cos{s},\sin{(u(t))}\sin{s}\right),
\end{aligned}
\end{equation}
where $t\in I\subset\r$ and $s\in\r$. We now calculate the mean curvature of $S_\gamma$.  The expression of the mean curvature  of $S_\gamma$ computed with the aid of \eqref{para1} is  
 \begin{equation}\label{mean}
 H= \frac{g_{11}h_{22}-2g_{12}h_{12}+g_{22} h_{11}}{2(g_{11}g_{22}-g_{12}^2)},
 \end{equation}
where, as  usual,   $\{g_{11},g_{12},g_{22}\}$ and $\{h_{11},h_{12},h_{22}\}$ are the coefficients of the first and second fundamental forms of $S_\gamma$ for the parametrization \eqref{para1}: 
 $$g_{11}=\langle \Phi_{t},\Phi_{t}\rangle,\quad g_{12}=\langle \Phi_{t},\Phi_{s}\rangle,\quad g_{22}=\langle \Phi_{s},\Phi_{s}\rangle,$$
 $$h_{11}=\langle G,\Phi_{tt}\rangle,\quad h_{12}=\langle G,\Phi_{ts}\rangle,\quad h_{22}=\langle G,\Phi_{ss}\rangle.$$
 Here $G$ is the unit normal vector field on $S_\gamma$. Note that $G(\Phi(t,s))$  is not only orthogonal to $\Phi_t(t,s)$ and $\Phi_s(t,s)$,  but also to $\Phi(t,s)$ since $G$ is a tangent vector of $\s^3$. A straightforward computation gives $g_{12}=h_{12}=0$,  
 $g_{11}=|\gamma'|^2$, $g_{22}=(\sin{u})^2$ and
 $$h_{11}=\frac{4u'^2\sin{u}+\cos{u}(2u''+\sin(2u))}{2|\gamma'|},\quad h_{22}=-\frac{\sin{u}(\cos{u})^2}{|\gamma'|}.$$
 Thus the mean curvature $H$ in \eqref{mean} is  
\begin{equation}\label{hs1}
 H=\frac{\cos{u}(\cos{u}+\cos{3u})+(3\cos{2u}-1)u'^2-2u''\sin{u}\cos{u} }{2 \sin{u}|\gamma'|^3}.
\end{equation}
The expression \eqref{ks} when $v(t)=t$ is    
\begin{equation}\label{kuu}
\kappa_s=\frac{1}{|\gamma'|^3}\left(2u'^2 \sin{u}+(\cos{u})^2\sin{u} +u''\cos{u}\right).
\end{equation}
Equation \eqref{kuu} allows to write $u''$ in terms of  $\kappa_s$. By replacing $u''$  in \eqref{hs1},  the mean curvature $H$ becomes
\begin{equation}\label{hh}
H=\frac{ (\cos{u})^2 -\sin{u}|\gamma'|\kappa_s }{ \sin{u}|\gamma'|}.
\end{equation}
Thus, $H=0$ if and only if $(\sin{u})|\gamma'|\kappa_s=(\cos{u})^2$. Therefore if $\gamma$ is a spherical catenary, the surface $S_\gamma$ is not minimal.

Looking in the formula \eqref{hh}, we observe that  the term $\sin{u}$ in  the numerator  is just  the Euclidean distance of the point $\Psi(u,v)$ to the plane $\Pi$ of equation $z=0$. This suggests to consider the potential energy of $\gamma$ calculated with respect to the plane $\Pi$ instead of the geodesic $P$. Definitively, we will formulate a different hanging chain problem in $\s^2$ in such a way that the critical points of the corresponding energy functional can successfully answer to the question of the minimality of $S_\gamma$.

Consider  a plane $\Pi$ of $\r^3$, which we can assume that it is the plane of equation $z=0$. As usual, the $z$-axis is  the direction of the gravity when the gravitational vector field is $\partial_z$.  Now, we replace the sphere $\s^2$ by an arbitrary surface $S$ of the Euclidean space $\r^3$. The {\it extrinsic hanging chain problem in $S$} consists in determining the shape of a hanging chain supported on $S$ where the potential energy of $\gamma$ is calculated with the Euclidean distance  to $\Pi$.  A critical point of this potential will be called an {\it extrinsic catenary on $S$}. Notice that the vector field $\partial_z$ is not a vector of $S$ but of the ambient space $\r^3$. In particular, coming back to the Euclidean context, now the chain in $S$ is subjected to the Euclidean gravity, which is constant. In particular,  we can assert that the potential at $ds$ of the chain is $\sigma g z\, ds$ as usual. 

The extrinsic hanging chain problem  was studied in the XIX century by Bobillier \cite{bo},  although it has not yet received much interest in the literature. See also \cite[Ch. VII]{ap} and \cite{gu}, and more recently,  \cite{fe}.  In the particular case that $S$ is the unit sphere $\s^2$, a solution of this problem will be called an {\it extrinsic spherical catenary} (Bobillier coined  the expression  ``spherical cha\^{i}nette'').  In this paper, we recall this problem and its solution  and, in addition, we credit to the work of Bobillier, an almost forgotten French mathematician  \cite{ha}. The potential energy of the hanging chain   is   
\begin{equation}\label{eex}
\mathcal{E}_S^{ex}[\gamma]=\int_a^b(\sin{u}+\lambda) |\gamma'(t)|\, dt= \int_a^b (\sin{u}+\lambda)\sqrt{u'^2+v'^2(\cos{u})^2}\, dt,
\end{equation}
where again $\lambda$ is a Lagrange multiplier. 

\begin{proposition} \label{pr1}
Let $\gamma$ be a regular curve in $\s^2_{+}$. Then  $\gamma$ is an extrinsic spherical catenary  if and only if its geodesic curvature $\kappa_s$ satisfies
\begin{equation}\label{ex-cat}
\kappa_s= \frac{v'(\cos{u})^2}{(\sin{u}+\lambda)|\gamma'|},
\end{equation}
or equivalently, if 
\begin{equation}\label{k-ex}
\kappa_s= \frac{\langle{\bf n},\partial_z\rangle}{\sin{u}+\lambda}.
\end{equation}
\end{proposition}
 
\begin{proof}
The Euler-Lagrange equations for \eqref{eex} follow   directly from \eqref{silva}, where now $f(u)=\sin{u}+\lambda$. This gives \eqref{ex-cat}. Formula \eqref{k-ex} is a consequence of  \eqref{ex-cat} and the expression \eqref{nn} for ${\bf n}$. 
\end{proof}

Equation \eqref{k-ex} is analogous to \eqref{catenary2} because the term  $\sin{u}$ in the denominator is the height with respect to $\Pi$ and the vector field $\partial_z$ is the gravitational vector field. 

Finally, we  answer the question of when the mean curvature of the rotational surface $S_\gamma$ is identically zero (objective T3).

\begin{corollary} \label{cor-s}
Let $\gamma$ be a regular curve in $\s^2_{+}\times\{0\}$.  Then $S_\gamma$ is minimal if and only if $\gamma$ is an extrinsic spherical catenary.
\end{corollary}
\begin{proof}
  Without loss of generality, we can write $\gamma$ as  $\gamma(t)=\Psi(u(t),t)$. From \eqref{hh}, the mean curvature $H$ vanishes if and only if $(\cos{u})^2=\sin{u}|\gamma'|\kappa_s$ and this identity is just  \eqref{ex-cat}.
\end{proof}
This result in $\s^3$ is analogous to the relation between the catenoid of $\r^3$ and the catenary curve obtained by Euler. 
Rotational surfaces in $\s^3$ with zero mean curvature (minimal surfaces) are known: see \cite{al,ri}. Among these surfaces, the Clifford torus is the most famous example because it is the only minimal embedded torus in $\s^3$ (\cite{br}). In the context of extrinsic spherical catenaries, the Clifford torus corresponds to the case $\kappa_s(t)= 1$ and $u(t)=\pi/4$ in \eqref{ex-cat}.  Indeed, the parametrization \eqref{para1} is 
$\Phi(t,s)=\frac{\sqrt{2}}{2}(\cos{t},\sin{t},\cos{s},\sin{s})$. Thus $S_\gamma=\s^1(\frac{1}{\sqrt{2}})\times \s^1(\frac{1}{\sqrt{2}})$ which it is the Clifford torus.

%%%%%%%%%%%%%%%%%%%%%%%%%%%%%%%
\section{The hanging chain problem in the hyperbolic plane}\label{sec3}
%%%%%%%%%%%%%%%%%%%%%%%%%%%%%%%%%%%%%%%%%%%%%%%%%%%%%%%

In this section,  the hanging chain problem in the hyperbolic plane $\h^2$ is investigated. The model for $\h^2$ will be the upper half-plane  $(\r^2_{+},g)$, where  $\r^2_+=\{(x,y)\in\r^2:y>0\}$ and the metric is  $g=\frac{dx^2+dy^2}{y^2}$. 

 The hanging chain problem in the hyperbolic plane is richer than in the Euclidean plane because there are several possibilities of  reference lines and potential energies. We will consider the situation that a horocycle is a   reference line. Horocycles have some analogies with the straight-lines of $\r^2$ and provide the so-called horospherical geometry  \cite{iz}.     This section is divided in three parts according to this variety of choices: 
\begin{enumerate}
\item Hyperbolic catenary: the reference line is a geodesic and the potential energy is calculated along geodesics of $\h^2$.
\item Hyperbolic horo-catenary: the reference line is a geodesic and the potential energy is calculated along horocycles of $\h^2$.
\item Horo-catenary: the reference line is a horocycle and the potential energy is calculated along geodesics of $\h^2$.
\end{enumerate}

Let $\gamma\colon [a,b]\to\h^2$ be a regular curve parametrized by $\gamma(t)=(u(t),v(t))$. The energy to minimize in all these situations in this section is of type 
\begin{equation}\label{ww}
\gamma\longmapsto   \int_a^b\omega(u,v)|\gamma'(t)|\, dt= \int_a^b\omega(u,v)\frac{\sqrt{u'^2+v'^2}}{v}\, dt
\end{equation}
where $\omega=\omega(u,v)$ is a smooth function on the variables $u$ and $v$. Here $ \sqrt{u'^2+v'^2}/v\, dt$ is the arc-length element of $\h^2$. This energy can be interpreted as the length of $\gamma$  in the conformal metric $\widetilde{g}= \omega^2g$ and consequently,   its critical points coincide with the  geodesics in the conformal space $(\r^2_{+},\widetilde{g})$. In order to simplify the presentation of this section,   the Euler-Lagrange equations of the energy \eqref{ww} are calculated in the following result.

\begin{proposition}\label{pr2} A regular curve  $\gamma(t)=(u(t),v(t))$ in $\h^2$  is a critical point of the energy \eqref{ww} if and only if its curvature $\kappa_h$ is
\begin{equation}\label{ayuda}
\kappa_h=\frac{v}{\omega\sqrt{m}}\left(u'\omega_v-v'\omega_u\right),
\end{equation}
where $m(t)=u'(t)^2+v'(t)^2$, $\omega_u=\frac{\partial\omega}{\partial u}$ and $\omega_v=\frac{\partial\omega}{\partial v}$. 
\end{proposition}

\begin{proof} A straight-forward computation of \eqref{sel} gives, respectively
$$\frac{u'}{v^2\sqrt{m}}\left(u'v\omega_v-v'\omega_u-\frac{u'\omega}{v}-\omega \frac{u'v''-v'u''}{m}\right)=0,$$
$$\frac{v'}{v^2\sqrt{m}}\left(u'v\omega_v-v'\omega_u-\frac{u'\omega}{v}-\omega \frac{u'v''-v'u''}{m}\right)=0.$$
Since $\gamma$ is regular, and using the Euclidean curvature $\kappa_e$ given in \eqref{k5}, we deduce that $\gamma$ is a critical point of the energy \eqref{ww} if and only if 
\begin{equation}\label{ayuda3}
\kappa_e=\frac{1}{\omega\sqrt{m}}\left(u'\omega_v-v'\omega_u-\frac{u'\omega}{v}\right).
\end{equation}
On the other hand, the curvature $\kappa_h$ of $\gamma$ is related to $\kappa_e$ because  the hyperbolic metric is conformal to the Euclidean one: see  \cite[Chapter 1]{be}.  This relation is 
\begin{equation}\label{keh}
\kappa_h= v\kappa_e+\frac{u'}{\sqrt{m}}.
\end{equation} 
Then \eqref{ayuda} is consequence of \eqref{ayuda3} and \eqref{keh}.
\end{proof}

The identity \eqref{ayuda} can be also expressed as follows.   Consider $\{\partial_x,\partial_y\}$ the canonical vector fields of $\r^2$. Then the gradient $\nabla\omega$  of $\omega$ (in $\h^2$) is 
$$\nabla\omega= y^2(\omega_x\partial_x+\omega_y\partial_y).$$ 
On the other hand, the unit normal vector of $\gamma$ is $\n(t)=v(t)\frac{(-v'(t),u'(t))}{\sqrt{m}}$. Thus we obtain:

\begin{corollary} A regular curve  $\gamma(t)=(u(t),v(t))$ in $\h^2$  is a critical point of the energy \eqref{ww} if and only if its curvature $\kappa_h$ satisfies
\begin{equation}\label{ayuda2}
\kappa_h=\frac{g(\n,\nabla\omega)}{\omega},
\end{equation}
where $\n$ is the unit normal vector of $\gamma$.  
  \end{corollary}
%%%%
\subsection{Hyperbolic catenaries}
%%%%%
The first case to investigate follows the same motivation as in the Euclidean plane. For the choice of the reference line, we take a geodesic $L$ of $\h^2$ which we can assume to be $L=\{(0,y):y>0\}$.  At this level, the potential will be $0$. The potential energy at each point  is determined by the hyperbolic distance to $L$ which is calculated along the geodesics orthogonal to $L$.   If  $(x,y)\in\h^2$, its distance  $d$ to $L$ is 
$$d=\log\frac{x+\sqrt{x^2+y^2}}{y}.$$
The geodesics orthogonal to $L$ are half-circles of $\r^2_{+}$ centered at the origin of $\r^2$. Thus the unit vector field $Y\in\mathfrak{X}(\h^2)$ which is orthogonal to all these geodesics at each point of $\h^2$ is    
\begin{equation}\label{yy}
Y(x,y)=y\left( \frac{ y}{\sqrt{x^2+y^2}}\, \partial_x-\frac{x}{\sqrt{x^2+y^2}}\, \partial_y\right).
\end{equation}
Given a curve $\gamma(t)=(u(t),v(t))$, define   the potential energy 
\begin{equation}\label{fh1}
\mathcal{E}_H[\gamma]=\int_a^b (d+\lambda)|\gamma'(t)|\, dt=\int_a^b (d+\lambda)\frac{\sqrt{u'^2+v'^2}}{v}\, dt, \quad  d=\log {\frac{u+r}{v}},
\end{equation}
where $r(t)=\sqrt{u(t)^2+v(t)^2}$.  As in the case of the sphere $\s^2$, in $\h^2$ we have no a notion of (constant) gravity. Let us observe that $\mathcal{E}_H$ is a particular case of \eqref{ww} by choosing $\omega(u,v)=d+\lambda$. It will be assumed that $d\not=0$, that is, $u\not=0$. Equivalently, the curve $\gamma$ is contained in one of the domains   $\h^2_{+}=\{(x,y)\in\h^2:x>0\}$ or $\h^2_{-}=\{(x,y)\in\h^2:x<0\}$. Since each domain is mapped into  other by means of the isometry  $(x,y)\mapsto (-x,y)$,   it will be assumed that $\mathcal{E}_H$  acts on the class of all curves $\gamma$ contained in $\h^2_{+}$.  

\begin{definition} A critical point of $\mathcal{E}_H$ is called a {\it hyperbolic catenary}.
\end{definition}

As in $\r^2$ and $\s^2$, a hyperbolic catenary will be characterized in terms of its curvature $\kappa_h$   as curve of $\h^2$.

\begin{theorem}\label{t32}
A regular curve  $\gamma(t)=(u(t),v(t))$ in $\h^2_+$  is a hyperbolic catenary if and only if its curvature $\kappa_h$ satisfies
\begin{equation}\label{eqh1}
\kappa_h=-\frac{1}{d+\lambda}\frac{uu'+vv'}{ \sqrt{u^2+v^2}\sqrt{u'^2+v'^2}}=-\frac{ r'}{(d+\lambda)\sqrt{u'^2+v'^2}}.
\end{equation}
\end{theorem}

\begin{proof}
The energy  $\mathcal{E}_H$ is a particular case of \eqref{ww}. Using \eqref{ayuda} with $f=d+\lambda$, then equation \eqref{eqh1} is obtained immediately.  
\end{proof}

With respect to T2, the next step consists of writing  equation \eqref{eqh1} in a similar manner as  the formula \eqref{catenary2} involving the curvature $\kappa_h$ and the vector field $Y$.  The following result is immediate by a direct computation or using \eqref{ayuda2} because the vector field $Y$ is just $\nabla d$.

\begin{corollary}\label{th1}
 A regular curve  $\gamma$  in $\h^2_{+}$   is a hyperbolic catenary if and only if its curvature $\kappa_h$ satisfies 
\begin{equation}\label{eth1}
\kappa_h= \frac{g( \n,Y)}{d+\lambda}.
\end{equation}
\end{corollary} 
As a consequence, Corollary \ref{th1} is the analogue in $\h^2$ of the statement   \eqref{catenary2} for hyperbolic catenaries.

%%%%%%%
\subsection{Hyperbolic horo-catenaries}
%%%%%%

Consider a modified version of the above  hanging chain problem replacing the potential calculated with the hyperbolic distance by the horocycle distance. The  {\it horocycle distance} to the geodesic $L$ is defined as the distance of a point $(x,y)\in\h^2$ to $L$ calculated by the horocyle  passing through  $(x,y)$ and orthogonal to $L$. In the present case that $L$ is the geodesic of equation $x=0$, this  distance   is  $|x|/y$. 

     Let  $\gamma(t)=(u(t),v(t))$, $t\in [a,b]$, be a regular curve contained in $\h^2_+$. The potential energy  of $\gamma$   calculated with  the horocycle distance is 
\begin{equation}\label{ehor}
\mathcal{E}_H^{hor}[\gamma]=\int_a^b (d_{hor}+\lambda) |\gamma'(t)|\, dt=\int_a^b (d_{hor}+\lambda) \frac{\sqrt{u'^2+v'^2}}{v}\, dt,\quad d_{hor}= \frac{u}{v}.
\end{equation}

\begin{definition} 
A critical point of $\mathcal{E}_H^{hor}$ is called an {\it hyperbolic horo-catenary}.
\end{definition}

With respect to  the objective T1, we prove:
 
\begin{theorem}
A regular curve $\gamma(t)=(u(t),v(t))$ in $\h^2_+$ is a  hyperbolic horo-catenary if and only if 
its curvature $\kappa_h$ satisfies 
\begin{equation}\label{h22}
\kappa_h=- \frac{uu'+vv'}{(d_{hor}+\lambda)v\sqrt{m}}.
\end{equation}
\end{theorem}

\begin{proof} The energy   $\mathcal{E}_H^{hor}$ in \eqref{ehor} is of type \eqref{ww} and formula \eqref{h22} is \eqref{ayuda} for $f= d_{hor}$. 
\end{proof}

To answer to T2, we replace the above vector field $Y$ by  the vector field $W\in\mathfrak{X}(\h^2)$ defined as 
$$W(x,y)=  y\, \partial_x- x\, \partial_y.$$
The next result is a consequence of \eqref{ayuda2} because $W=\nabla d_{hor}$.

\begin{corollary} A regular curve $\gamma$  in $\h^2_{+}$ is a hyperbolic horo-catenary   if and only if its curvature $\kappa_h$ satisfies
\begin{equation}\label{eth2}
\kappa_h= \frac{g( \n,W)}{d_{hor}+\lambda}.
\end{equation}
\end{corollary}

To conclude this subsection we investigate problem T3 for this type of catenaries. The hyperbolic plane $\h^2$ is embedded into the $3$-dimensional hyperbolic space $\h^3=(\r^3_{+},\frac{1}{x_3^2}(dx_1^2+dx_2^2+dx_3^2))$ via the natural inclusion  $(x,y)\in\h^2\mapsto (x ,0,y)\in\h^3$. With this identification, the geodesic $L\subset\h^2$ is the $x_3$-axis in $\h^3$.    
Let $S_\gamma$ denote the surface of revolution obtained by rotating $\gamma(t)=(u(t),0,v(t))$  with respect to the $x_3$-axis. In the upper half-space model of $\h^3$, the rotations that leave pointwise fixed the $x_3$-axis coincide with the Euclidean rotations of $\r^3$ with the same axis. These surfaces of revolution in $\h^3$  are called of spherical type  \cite{dcarmo}. Thus a parametrization $\Phi$  of $S_\gamma$ is 
\begin{equation}\label{para2}
\Phi(t,s)= (u(t)\cos{s},u(t)\sin{s},v(t)) , \quad t\in [a,b],s\in\r.
\end{equation}

\begin{theorem}\label{t37} A regular curve  $\gamma$   in $\h^2_+$   is a hyperbolic horo-catenary for $\lambda=0$ if and only the rotational surface $S_\gamma$ of spherical type is minimal.
\end{theorem}

\begin{proof}  In the upper half-space model of $\h^3$, the mean curvature $H$ of a surface $S$ can be computed with the aid of  the Euclidean  mean curvature $H_e$ of $S$  when $S$ is viewed as a submanifold of the Euclidean space $\r^3_{+}$. This relation is similar to \eqref{keh}, namely,  
\begin{equation}\label{he0}
H(p)=x_3H_e(p)+N_3(p),
\end{equation}
 where $p=(x_1,x_2,x_3)\in S$ and $N=(N_1,N_2,N_3)$ is the Euclidean unit normal vector of $S$ (\cite[Chapter 1]{be}). If now $S$ is the rotational surface $S_\gamma$ parametrized by $\Phi$ in  \eqref{para2}, the value of $H_e$  is 
\begin{equation}\label{he}H_e=\frac{\kappa_e}{2}+\frac{v'}{2u\sqrt{m}},
\end{equation}
and the expression of $N$ is   
$$N=\frac{1}{\sqrt{m}}(-v'\cos s,-v'\sin s,u').$$
Thus $N_3=u'/\sqrt{m}$, and using \eqref{he}, the mean curvature $H$ given in \eqref{he0} becomes 
$$H=\frac{v\kappa_e}{2}+\frac{vv'}{2u\sqrt{m}}+\frac{u'}{\sqrt{m}}.$$
Using \eqref{keh}, 
$$H=\frac{\kappa_h}{2}+\frac{uu'+vv'}{2\sqrt{m}}.$$
Then $H=0$ if and only if 
\begin{equation}\label{vke}
\kappa_h=-\frac{uu'+vv'}{\sqrt{m}}.
\end{equation}
But this identity \eqref{vke} is just   equation \eqref{h22} for $\lambda=0$   because $d_{hor}=u/v$. This proves the result. 
\end{proof}
 
 We point out that  do Carmo and Dajczer obtained   all  minimal rotational surfaces of $\h^3$.  The statement of Theorem \ref{t37} gives a geometric interpretation of the generating curves of minimal rotational surfaces of spherical type of $\h^3$ proving that these curves are the solutions of a hanging chain problem in $\h^2$.  As a consequence, this extends the Euler's result to spherical minimal rotational surfaces.
 %%%%%%%%%%%%%%
 \subsection{Horo-catenaries}
 %%%%%%%%%%%%%%%%%%%%%%%%%%%%%%%%%%%%%%%%%%%%%%%%%%%%%%
 We investigate the hanging chain problem  considering a   horocycle $\mathcal{H}$   as reference line.   Without loss of generality we can assume   $\mathcal{H}=\{(t,1):t\in\r\}$.  The potential energy at each point of $\h^2$ is given by its hyperbolic distance to $\mathcal{H}$. In the upper half-plane model of $\h^2$, the geodesics orthogonal to $\mathcal{H}$ are vertical lines of $\r_+^2$. If $(x,y)\in\h^2$, the hyperbolic distance $d_b$ from $(x,y)$ to $\mathcal{H}$ is  the length throughout the geodesic orthogonal to $\mathcal{H}$ passing through $(x,y)$. This distance is   $d_b=\log (y)$.  Note that this distance coincides with the Busemann function   in the horospherical geometry when the ideal point is $\infty$ \cite{bu}. The unit vector field  $V\in\mathfrak{X}(\h^2)$ which is tangent to all these geodesics  is given by 
  $$V(x,y)=y\partial_y.$$
       We will assume again that $d_b\not=0$, that is, $y\not=1$. The horocycle $\mathcal{H}$ separates $\h^2$ in two   domains, namely, $\h^2(+)=\{(x,y)\in\h^2:y>1\}$ and $\h^2(-)=\{(x,y)\in\h^2:y<1\}$, but both domains are not isometric. From now on, we will assume that all curves are contained in $\h^2(+)$ and a similar work can be done in the case that all curves are contained in $\h^2(-)$. 

Let $\gamma\colon [a,b]\to\h^2_+$ be a regular curve, $\gamma(t)=(u(t),v(t))$. The potential energy   of $\gamma$ is 
\begin{equation}\label{ehor2}
 \mathcal{E}_{hor}[\gamma]=\int_a^b  ( d_b+\lambda)|\gamma'(t)|\, dt=\int_a^b  ( d_b+\lambda) \frac{\sqrt{u'^2+v'^2}}{v}\, dt.\quad d_b=\log{v},
 \end{equation}
where   $\lambda\in\r$ is a Lagrange parameter. 

\begin{definition}
A critical point of $\mathcal{E}_{hor}$ is called a {\it horo-catenary}.
\end{definition}

We characterize the horo-catenaries in terms of their curvatures $\kappa_h$.

\begin{theorem} A regular curve $\gamma(t)=(u(t),v(t))$ in $\h^2(+)$ is a horo-catenary if and only if its curvature $\kappa_h$ satisfies
\begin{equation}\label{h33}
\kappa_h= \frac{u'}{(d_b+\lambda)\sqrt{m}}.
\end{equation}
\end{theorem}

\begin{proof}
Expression \eqref{h33} is just \eqref{ayuda} for $f=d_b+\lambda$. 
\end{proof}

The Lagrangian    $J$ of $\mathcal{E}_{hor}$ is 
\begin{equation}\label{juv}
J[u,v,u',v']=( d_b+\lambda) \frac{\sqrt{u'^2+v'^2}}{v}, 
\end{equation}
which does not depend on $u$. Thus  a first integration of the Euler-Lagrange equation can be deduced.

\begin{corollary}  A regular curve $\gamma$  in $\h^2(+)$  is a horo-catenary if and only if $\gamma$ can be locally expressed as
\begin{equation}\label{integra1}
\gamma(v)=\left(\int^v\frac{c\tau}{\sqrt{(\log{\tau}+\lambda)^{2}-c^2 \tau^2}}\, d\tau,v\right),
\end{equation} 
where $c\in\r$ is a constant of integration. 
\end{corollary}
\begin{proof}  From   \eqref{juv}, there exists a constant $c$ such that  $\frac{\partial J}{\partial u'}=c$. This identity is 
\begin{equation}\label{cc}
\frac{u'(\log{v}+\lambda)}{v\sqrt{u'^2+v'^2}}=c.
\end{equation}
Without loss of generality, we can assume that $\gamma$ writes locally as $\gamma(v)=(u(v),v)$. Then \eqref{cc} is 
$$\frac{u'(\log{v}+\lambda)}{v\sqrt{1+u'^2}}=c.$$
Hence it follows \eqref{integra1}.\end{proof}

With respect to the objective T2, we have the following characterization of horo-catenaries which is a consequence of \eqref{ayuda2} and the fact that $V=\nabla d_b$.

\begin{corollary} A regular curve $\gamma$ in $\h^2(+)$ is a horo-catenary if and only if  its curvature $\kappa_h$ satisfies
\begin{equation}\label{kh2}
\kappa_h = \frac{g( {\bf n},V)}{d_b+\lambda}.
\end{equation}
\end{corollary}

Equation \eqref{kh2} is the analogue of  the formula \eqref{catenary2} in the context of horo-catenaries.

\begin{figure}[hbtp]
\centering
\includegraphics[width=.3\textwidth]{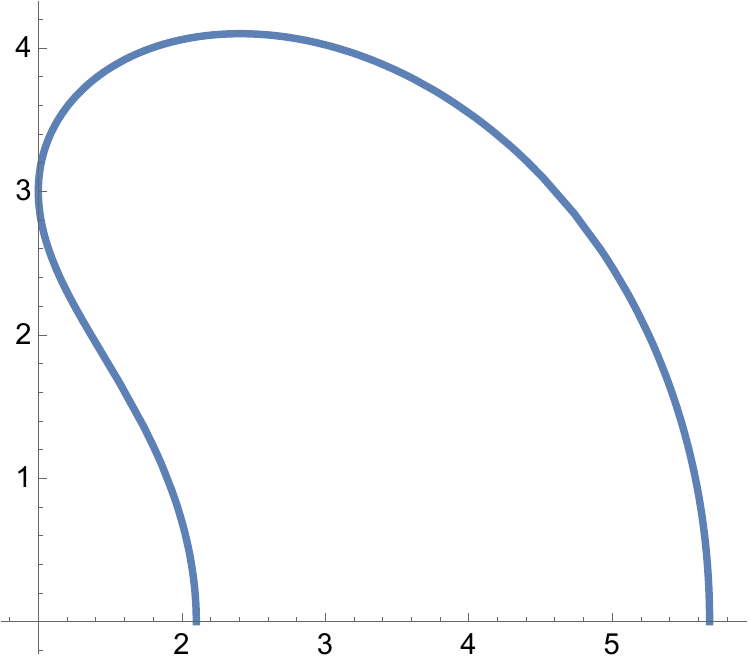} \includegraphics[width=.25\textwidth]{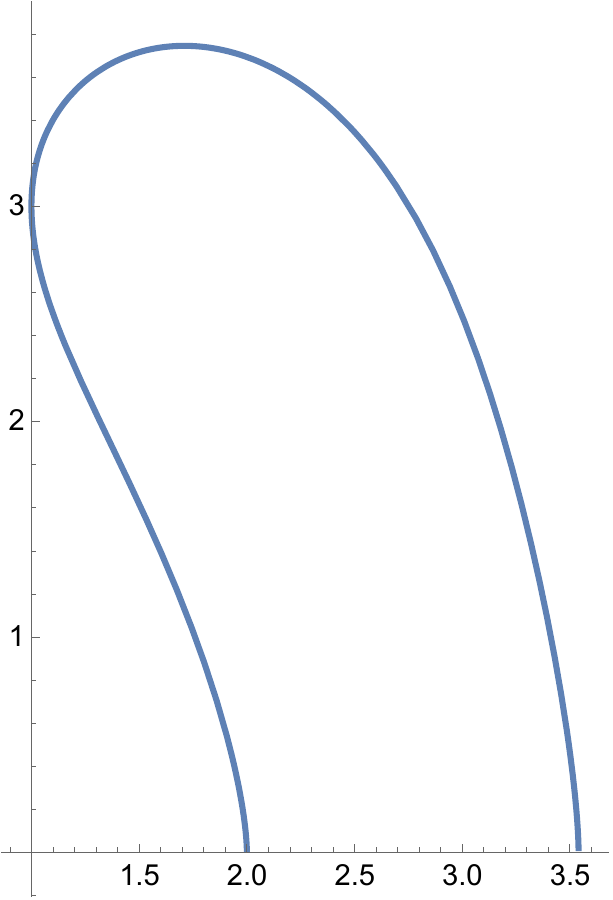} \includegraphics[width=.4\textwidth]{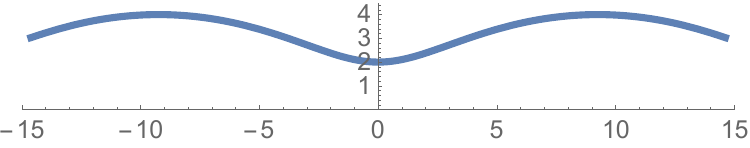}
\caption{Catenaries in $\h^2$ considering the upper half-plane model: hyperbolic catenary (left), hyperbolic horo-catenary (middle), horo-catenary (right). These curves are solutions of \eqref{ode} and the initial conditions \eqref{conditions} are: 
$u_0=1$, $v_0=3$ and $\theta_0=\pi/2$ (left and middle); $u_0=0$, $v_0=1.8$ and $\theta_0=0$ (right).} \label{fig2}
\end{figure}

We finish this section   with some pictures of the three types of catenaries  for $\lambda=0$. See Figure \ref{fig2}. The process to plot these curves with {\it Mathematica} is the following. Suppose that $\gamma(t)=(u(t),v(t))$ is parametrized by the Euclidean arc-length. Then $\gamma'(t)= (\cos\theta(t),\sin\theta(t))$ for some function $\theta=\theta(t)$, where $\theta'(t)=\kappa_e(t)$ is the Euclidean curvature of $\gamma$. For each type of catenary, the value of $\kappa_e$ is obtained by combining \eqref{keh} and each one of the expressions for $\kappa_h$ in  \eqref{eqh1}, \eqref{h22} and \eqref{h33}:  
\begin{align}
\kappa_e&=-\dfrac{1}{v\sqrt{m}}\left( \dfrac{uu'+vv'}{rd}+u'\right),&\quad (\mbox{hyperbolic catenary}),\label{h11}\\
\kappa_e&=-\dfrac{1}{v\sqrt{m}}\left( \dfrac{uu'+vv'}{vd_{hor}}+u'\right),&\quad (\mbox{hyperbolic horo-catenary}),\label{h3}\\
\kappa_e&=\dfrac{1}{v\sqrt{m}}\dfrac{1-d_b}{d_b}u',&\quad (\mbox{horo-catenary})\label{h-33}.
\end{align}
 
 It follows that the functions $u(t)$, $v(t)$ and $\theta(t)$ satisfy the ODE system
\begin{equation}\label{ode}
\left\{ \begin{aligned}
u'(t)&=\cos\theta(t)\\
v'(t)&=\sin\theta(t)\\
\theta'(t)&=\kappa_e(t).
\end{aligned}
\right.
\end{equation}
Finally, and distinguishing the three types of catenaries of $\h^2$, the system \eqref{ode} have been numerically solved with {\it Mathematica} once initial conditions 
\begin{equation}\label{conditions}
u(0)=u_0,\quad v(0)=v_0,\quad \theta(0)=\theta_0
\end{equation}
 have been prescribed.

As a consequence of these plots, we observe that the horo-catenary (Figure \ref{fig2}, right)  is a graph on the $x$-axis. This is not   a coincidence, but it holds in general.  

\begin{proposition}\label{pr-h}
 If $\gamma$ is a horo-catenary, then $\gamma$ is a vertical line or $\gamma$ is a  bounded entire graph on the $x$-axis.
\end{proposition}

\begin{proof}
Suppose that $\gamma$ is parametrized by $\gamma(t)=(u(t),v(t))$, $t\in I\subset\r$, where $t$ is the Euclidean arc-parameter and   $I$ is the maximal domain.  Since $\gamma$ is a horo-catenary, the curvature $\kappa_e$ of $\gamma$ satisfies \eqref{h-33}. Thus the ODE system \eqref{ode} is  
\begin{equation}\label{s-hor}
\left\{ \begin{aligned}
u'(t)&=\cos\theta(t)\\
v'(t)&=\sin\theta(t)\\
\theta'(t)&=\frac{1-\log{v}}{v\log{v}}u'.
\end{aligned}
\right.
\end{equation}
We distinguish two cases.
\begin{enumerate}
\item Suppose there exists $t_0$ such that $u'(t_0)=0$.  The first equation of \eqref{s-hor} implies $\cos\theta(t_0)=0$. Without loss of generality, we suppose that $\theta(t_0)=\pi/2$. Thus at $t=t_0$, the initial conditions \eqref{conditions} are $(u(t_0),v(t_0),\pi/2)$. By uniqueness of \eqref{s-hor}-\eqref{conditions}, $\gamma$ is a vertical straight-line.  
\item Suppose $u'(t)\not=0$ for all $t$. This implies that $\gamma$ is a graph on some interval $I=(a,b)$ of the $x$-axis. Reparametrizing $\gamma$,  the curve $\gamma$ can be expressed by $\gamma(x)=(x,v(x))$, $x\in I$ and it will be proved that $I=\r$. Equation \eqref{cc} for $\kappa_e$ becomes 
\begin{equation}\label{log}
\frac{\log{v}}{v\sqrt{1+v'^2}}=c.
\end{equation}
Since $c\not=0$ and $\gamma$ is contained in $\h^2(+)$, we have $\log(v)>0$. From \eqref{log} we deduce
\begin{equation}\label{logv}
0<c=\frac{\log{v}}{v\sqrt{1+v'^2}}\leq \frac{\log{v}}{v}.
\end{equation}
The function $t\mapsto \log{t}/t $ is bounded in $(1,\infty)$ with the property 
\begin{equation}\label{log2}
\lim_{t\to\infty} \frac{\log{t}}{t}=\lim_{t\to 1} \frac{\log{t}}{t}=0.
\end{equation}
From \eqref{logv} and \eqref{log2}, we conclude that $v'$ is a bounded function. Moreover, there exist two constants $m_1, m_2\in\r$ such that $1<m_1<m_2$ and $m_1\leq v(x)\leq m_2$ for all $x\in (a,b)$. The fact that the function $v'(t)$ is bounded proves finally that all solutions of \eqref{s-hor} are defined in the entire real line $\r$. 
\end{enumerate}
\end{proof}

%%%%
\section{Conclusions and outlook}

The catenary is the solution of the hanging chain problem in $\r^2$ and this makes it so attractive in other fields of science, engineering and  architecture. However,   the hanging chain problem has not been formulated in spaces other than Euclidean one. Among these spaces,  the sphere $\s^2$  and the hyperbolic plane $\h^2$ are the natural choices to extend this problem. It has been formulated this problem in $\s^2$ and in $\h^2$, defining in each case a potential energy that depends on the distance of a point with respect to a reference line. The resulting critical points of these energies (for different reference lines) have generalized the concept of catenary in both spaces.

A remarkable result is the characterization of the generating curves of minimal rotational surfaces of $\s^3$ proving that these curves are chains on $\s^2$ suspended by its weight where the force vector field is really the gravity of $\r^3$. In this particular situation, the initial hanging chain problem formulated in $\s^2$ must be replaced by other, which was called `extrinsic', because the force field  is a vector of $\r^3$, not of $\s^2$. Then  it was proved that the generating curves are solutions of an old problem formulated by Bobillier in the 19th century  and that it has been revisited in the present paper.

There are a number of problems in which this article could be expanded. For example, a question concerns to investigate the existence of closed spherical catenaries. In view of the pictures of Figure \ref{fig1}, it seems plausible that such catenaries do exist. This problem was investigated in \cite{al} for extrinsic spherical catenaries in the context of rotational surfaces of $\s^3$ with constant mean curvature.  Besides the closed catenaries,  there are other catenaries which never close up  and they are turning around the north pole of $\s^2$. The problem that arises here is that, although the curvature function is periodic, this is not enough to ensure that the corresponding catenary will be closed: see a discussion of this problem in \cite{arroyo}.

As in $\s^2$,  it would be interesting to classify the catenaries in the hyperbolic plane. Proposition \ref{pr-h} is just an example, but the work to be done goes beyond that. According to Figure \ref{fig2}, several questions are reasonable to ask. For example, (i) when does a catenary intersect the ideal boundary of $\h^2$? and in  such a case, determine whether the intersection is orthogonal; (ii) is  every  horo-catenary periodic? (iii) which are the properties of the  horo-catenaries contained in $\h^2(-)$?

Another extension of the paper would be to consider the shape of a hanging surface in $\s^3$ and $\h^3$.   In the Euclidean space, the analogue of the catenary in the two-dimensional case is called a singular minimal surface  (\cite{bht,dih}). The extension is straightforward using the characterization \eqref{catenary2}. So, it suffices to replace the curvature of the catenary $\kappa_e$ by the mean curvature $H$ of the surface and  the unit normal ${\bf n}$ of the curve by the unit normal vector field $G$ to the surface. For example, in the $3$-dimensional unit sphere   $\s^3$, the shape of a hanging surface  with respect to $\s^2\times\{0\}$ is characterized by the equation
$H= \langle G,X\rangle/(d+\lambda)$, $\lambda\in\r$. Here $d$ is the distance to $\s^2\times\{0\}$ and $X\in\mathfrak{X}(\s^3$) is the unit vector field tangent to the meridians of $\s^3$ which are orthogonal to $\s^2\times\{0\}$.

Finally, it could be interesting to obtain some geometric properties of the rotational surfaces in $\s^3$ and in $\h^3$ constructed by catenaries in its different possibilities. Although, initially the hanging chain problem has no a relation to the problem for rotational surfaces with minimum area, in some cases we have proved a connection between both problems (Corollary \ref{cor-s} and  Theorem \ref{t37}). It seems interesting to investigate geometric properties of the rotational surfaces of $\s^3$ and $\h^3$ whose generating curves are catenaries of $\s^2$ and $\h^2$, respectively.  

\section*{Declaration of competing interest} The author declares that he has no known competing financial interests or personal relationships that could have appeared to influence the work reported in this paper.

\section*{Acknowledgements}  The    author  is a member of the Institute of Mathematics  of the University of Granada. This work  has been partially supported by  the Project  PID2020-117868GB-I00 and MCIN/AEI/10.13039/501100011033.
 %%%%%%%%%%%%%%%%%%%%%%%%%%%%%%%%%%%%%%
 %%%%%%%%%%%%%%%%%%%%%%%%%%%%%%%%%%%%%%
 
 %%%%%%%%%%%%%%%%

 \end{document}